\documentclass[12pt]{article}
\usepackage[top=1in,bottom=1in,left=1in,right=1in]{geometry}
\usepackage{indentfirst}
\usepackage{amsfonts,amsmath,amsthm,amssymb,hyperref,enumerate,bm}
\usepackage{longtable}
\usepackage[caption=false]{subfig}
\usepackage{leftidx}
\usepackage{lineno}
\usepackage{color,xcolor}
\usepackage{array}
\usepackage{latexsym}
\usepackage{float}
\usepackage{color}
\usepackage{longtable,supertabular,booktabs}
\usepackage{rotating}
\usepackage{cases}
\usepackage{multicol,multirow}
\usepackage{epsfig}
\usepackage{cite}
\usepackage{fancyhdr}       
\usepackage{dsfont}
\usepackage{hyperref}
\hypersetup{
	colorlinks=true,
	linkcolor=red,
	filecolor=blue,
	urlcolor=red,
	citecolor=blue,
}

\newtheorem{theorem}{Theorem}[section]

\newtheorem{lemma}{Lemma}[section]

\newtheorem{corollary}{Corollary}[section]
\newtheorem{remark}{Remark}[section]

\def\G1{G^\mathcal{C}}

 \newenvironment{prof}{\trivlist
      \item[\hskip\labelsep
      {\itshape Proof.}]\normalfont}
      {\hspace*{\fill}$\Box$\endtrivlist}

\begin{document}
\title{The order of appearance of the product of the first and second Lucas numbers}
\author{
Hongjian Li
\footnote{E-mail\,$:$ lhj@gdufs.edu.cn. Supported by the Project of Guangdong University of Foreign Studies (Grant No. 2024RC063).}\\
{\small\it  School of Mathematics and Statistics, Guangdong University of Foreign Studies,}\\
{\small\it Guangzhou 510006, Guangdong, P. R. China} \\
Huiming Xiao
\footnote{Corresponding author. E-mail\,$:$ huiming\_xiao@m.scnu.edu.cn.}\quad
Pingzhi Yuan
\footnote{E-mail\,$:$ yuanpz@scnu.edu.cn. Supported by the National Natural Science Foundation of China (Grant No. 12171163) and the Basic and Applied Basic Research Foundation of Guangdong Province (Grant No. 2024A1515010589).}\\
{\small\it  School of Mathematical Sciences, South China Normal University,}\\
{\small\it Guangzhou 510631, Guangdong, P. R. China} \\
}
\date{}
 \maketitle
\date{}

\noindent{\bf Abstract}\quad
Let $a$ and $b$ be relatively prime integers. Then the first Lucas sequence $\left(U_n\right)_{n\geq0}$ and the second Lucas sequence $\left(V_n\right)_{n\geq0}$ are defined respectively by $U_{n+2}=aU_{n+1}+bU_{n},\, U_0=0,\,U_1=1$ and $V_{n+2}=aV_{n+1}+bV_{n},\, V_0=2,\,V_1=a$, where $n\geq0$. Let $m$ be an integer with $\gcd(m,\,b)=1$. Then the smallest positive integer $k$ satisfying $m\mid U_k$ is called the order of appearance of $m$ in the first Lucas sequence $(U_n)_{n\geq0}$, denoted by $\tau(m)$, i.e., $\tau(m):=\min\{k\geq1:m\mid U_k\}$. When $a>0$ and $\Delta=a^2+4b>0$, we give explicit formulae for $\tau(U_m V_n), \tau(U_m U_n)$, $\tau(V_m V_n)$ and $\tau(U_nU_{n+p}U_{n+2p})$, thus generalizing the results of Irmak and Ray \cite{Irmak}.

%
%


\medskip
\noindent{\bf MR(2020) Subject Classification} 11B39, 11A05

\section{Introduction}\label{sec1}
Let $a$ and $b$ be relatively prime integers. Then for each integer $n\geq0$, define $U_n=U_n(a,b)$ and $V_n=V_n(a,b)$ as follows:
\begin{equation}\label{eq7.15.1}
U_{n+2}=aU_{n+1}+bU_{n},\quad U_0=0,\,U_1=1
\end{equation}
and
\begin{equation}\label{eq7.15.2}
V_{n+2}=aV_{n+1}+bV_{n},\quad V_0=2,\,V_1=a.
\end{equation}
The sequences $U=(U_n(a,b))_{n\geq0}$ and $V=(V_n(a,b))_{n\geq0}$ are called the first and second Lucas sequences with parameters $(a,b)$, respectively\cite{Ribenboim}.
If we take $a=b=1$, then the numbers $U_n=U_n(1,1)$ are called the Fibonacci numbers, while the numbers $V_n=V_n(1,1)$ are called the Lucas numbers.
Let $\alpha$ and $\beta$ be the roots of the characteristic equation $x^2-ax-b=0$, i.e., $\alpha=\left(a+\sqrt{\Delta}\right)/2$ and $\beta=\left(a-\sqrt{\Delta}\right)/2$, where $\Delta=a^2+4b$ is the discriminant. Then the Binet's formulae for the first and second Lucas sequences are
\begin{equation}\label{eq1}
U_n=\frac{\alpha^n-\beta^n}{\alpha-\beta}
\end{equation}
and
\begin{equation}\label{eq2}
V_n=\alpha^n+\beta^n,
\end{equation}
respectively.
For convenience, we mainly consider the non-degenerate first and second Lucas sequences.
That is, $b\neq0$ and the ratio $\alpha/\beta$ is not a root of unity,  which excludes the pairs $(a,b)\in \{(\pm2,-1), (\pm1, -1), (0,\pm1), (\pm1, 0)\}$ \cite[pp. 5-6]{Ribenboim}.
Throughout this paper, we let $(U_n)_{n\geq0}$ and $(V_n)_{n\geq0}$ denote the non-degenerate first and second Lucas sequences  with parameters $(a,b)$, respectively.

Let $m$ be a positive integer, which is relatively prime to $b$. Then the order of appearance of $m$ in the first Lucas sequence $(U_n)_{n\geq0}$ is defined as the smallest positive integer $k$ such that $m$ divides $U_k$ and denoted by $\tau(m)$ \cite{Renault}. If we take $a=b=1$, then $\tau(m)$ is the order of appearance of $m$ in Fibonacci sequence $(F_n)_{n\geq0}$, and we denote it by $z(m)$ in this case. There are many results about $\tau(m)$ and $z(m)$ in the literature.  Marques \cite{5Marques,6Marques,7Marques,8Marques} established explicit formulae for
\[
z(F_m\pm1),\,z(F_nF_{n+1}F_{n+2}F_{n+3}),\, z(F_n^{k+1})\mbox{~and~}z(L_{n}L_{n+1}L_{n+2}L_{n+3}),
\]
where $(L_n)_{n\geq0}$ is the classical Lucas sequence.
Pongsriiam \cite{Pongsriiam} obtained the formula for $z(L_n^k)$ for all $n,k\geq 1$.
Later,  Khaochim \cite{Khaochim1, Khaochim2} derived new results for $z(F_nF_{n+1}\ldots F_{n+k})$
and $z(L_nL_{n+1}\ldots L_{n+s})$, where $k \in \{4,5,6\}$ and  $s \in \{5,6\}$.
Trojovsk\'{y} \cite{1Trojovsky, 2Trojovsky}  gave explicit formulae for $z(L_m-L_n)$ and $z(F_m\pm F_n)$.
Ray et al.\cite{Ray} derived an explicit formula for $\tau(U_n^{k+1})$ with $k \geq 0$.
Irmak and Ray  \cite{Irmak} derived explicit formulae for
\[
z(F_m L_n),\,z(F_m F_n),\,z(L_m L_n)\mbox{~and~}z(F_nF_{n+p}F_{n+2p}),
\]
where $p\geq 3$ is a prime number. In this paper, we generalize the results of Irmak and Ray and give the explicit formulae for
\[
\tau(U_m V_n),\, \tau(U_m U_n),\, \tau(V_m V_n) \mbox{~and~} \tau(U_nU_{n+p}U_{n+2p}),
 \]
 where  $p\geq3$ is a  prime number. That is, we have the following theorems.

\begin{theorem}\label{Main theorem}
Let $(U_n)_{n\geq0}$ and $(V_n)_{n\geq0}$ be the first and second Lucas sequences as in \eqref{eq7.15.1} and \eqref{eq7.15.2}, respectively.
Let $m,\,n\geq 3$ be integers. Let $d=\gcd(m,\,n)$ and let $[m,n]$ denote the least common multiple of the integers $m$ and $n$.
Furthermore, let $\nu_p(n)$ denote the $p$-adic valuation of $n$.
If $a>0$ and $\Delta=a^2+4b>0$, then we have the following statements.
\begin{enumerate}
\item[\rm(i)]\, $\tau(U_mV_n)=\left\{
\begin{aligned}
&2[m,n],\  &&\text{if} \  \nu _2(m)\leq \nu_2(n),\\
&[m,n]V_d,\         &&\text{if} \  \nu_2(m)>\nu_2(n);
\end{aligned}
                   \right.
                     $
\item[\rm(ii)]\, $\tau(U_mU_n)=[m,n]U_d$;
\item[\rm(iii)]\,
{\footnotesize
$
\tau(V_mV_n)=\left\{
\begin{aligned}
&[m,n]\gcd(V_m,\,V_n),\    && \text{if}\ 2\nmid b,\,2\mid a,\,2\mid d\ \text{or}\ 2\nmid b,\,2\nmid a,\,3\mid d\ \text{or}\ 2\nmid b,\,2\mid a,\,2\nmid d,\,\nu_2(m)\neq\nu_2(n),\\
&2[m,n]\gcd(V_m,\,V_n), &&\text{if}\ 2\mid b\ \text{or}\ 2\nmid b,\,2\nmid a,\,3\nmid d\ \text{or}\ 2\nmid b,\,2\mid a,\,2\nmid d,\,\nu_2(m)=\nu_2(n).
\end{aligned}
                   \right.
$
}
\end{enumerate}

\end{theorem}

Let $a=b=1$ in Theorem \ref{Main theorem}. Then we have the following corollary, which is given by \cite[Theorem 1.1]{Irmak}.

\begin{corollary}\label{corollary 1}
Let $(F_n)_{n\geq0}$ and $(L_n)_{n\geq0}$ be the classical Fibonacci and Lucas sequences, respectively. Let $m,\,n \geq 3$ be integers and let $d=\gcd(m,\,n)$. Then the following statements hold.
\begin{enumerate}
\item[\rm(i)]\, $z(F_mL_n)=\left\{
\begin{aligned}
&2[m,n],\  &&\text{if}\ \nu_2(m)\leq \nu_2(n),\\
&[m,n]L_d,\         &&\text{if}\ \nu_2(m)>\nu_2(n);
\end{aligned}
                   \right.
                     $
\item[\rm(ii)]\, $z(F_mF_n)=[m,n]F_d$;
\item[\rm(iii)]\, $z(L_mL_n)=\left\{
\begin{aligned}
&[m,n]\gcd (L_m,\,L_n),\  &&m\equiv n\equiv 0\pmod3,\\
&2[m,n]\gcd (L_m,\,L_n),\         &&\text{otherwise}.
\end{aligned}
                   \right.
                     $
\end{enumerate}
\end{corollary}

\begin{theorem}\label{Th1.2}
Let $(U_n)_{n\geq0}$ and $(V_n)_{n\geq0}$ be the first and second Lucas sequences as in \eqref{eq7.15.1} and \eqref{eq7.15.2}, respectively.
If $a>0$ and $\Delta=a^2+4b>0$, then for each positive integer $n$ and each prime number $p\geq3$, we have
\[
\tau(U_nU_{n+p}U_{n+2p})=\left\{
\begin{aligned}
&n(n+p)(n+2p), &\text{if}~&\  p\nmid n\    \text{and}\  2\nmid n, \\
&\frac{n(n+p)(n+2p)}{2}\cdot\frac{a}{\gcd\left(a,\,n+p\right)}, &\text{if}~&\  p\nmid n\   \text{and}\  2\mid n, \\
&\frac{n(n+p)(n+2p)}{p^2}U_p^2, &\text{if}~&\  p\mid n\    \text{and}\  2\nmid n, \\
&\frac{n(n+p)(n+2p)}{2p^2}\cdot\frac{U_p^2V_p}{\gcd\left(V_p,\, (n+p)/p\right)}, &\text{if}~&\  p\mid n\    \text{and}\  2\mid n.
\end{aligned}
\right.
\]
\end{theorem}
Let $a=b=1$ in Theorem \ref{Th1.2}. Then we have the following corollary, which is given by \cite[Theorem 1.2]{Irmak}.

\begin{corollary}\label{corollary 2}
Let $(F_n)_{n\geq0}$ and $(L_n)_{n\geq0}$ be the classical Fibonacci and Lucas sequences, respectively.
Then for each positive integer $n$ and each prime number $p\geq3$, we have
\[
z(F_nF_{n+p}F_{n+2p})=\left\{
\begin{aligned}
&n(n+p)(n+2p), &\text{if}~&\  p\nmid n\    \text{and}\  2\nmid n, \\
&\frac{n(n+p)(n+2p)}{2}, &\text{if}~&\  p\nmid n\   \text{and}\  2\mid n, \\
&\frac{n(n+p)(n+2p)}{p^2}F_p^2, &\text{if}~&\  p\mid n\    \text{and}\  2\nmid n, \\
&\frac{n(n+p)(n+2p)}{2p^2}\cdot\frac{F_p^2L_p}{\gcd\left(L_p,\, (n+p)/p\right)}, &\text{if}~&\  p\mid n\    \text{and}\  2\mid n.
\end{aligned}
\right.
\]
\end{corollary}
\begin{remark}
\rm
The result  of Irmak and Ray\cite[Theorem 1.2]{Irmak} is different from our result of Corollary \ref{corollary 2} in the case when $p\mid n$ and $2\mid n$.
A counterexample is presented below.
If $n=50$ and $p=5$, then $F_5=5$, $L_5=11$ and
\begin{align}
z(F_nF_{n+p}F_{n+2p}) &=\frac{n(n+p)(n+2p)}{2p^2}\cdot\frac{F_p^2L_p}{\gcd\left(L_p,\, (n+p)/p\right)}  \notag\\
                      &=\frac{50\cdot55\cdot60}{2\cdot5^2}\cdot\frac{5^2\cdot11}{\gcd(11,\ 11)}\notag\\
                   &=25\cdot55\cdot60 \notag\\
                   &<25\cdot55\cdot60\cdot11  \notag\\
                   &=\frac{n(n+p)(n+2p)}{2p^2}F_pF_{2p}. \notag
\notag \end{align}
\end{remark}

\section{Preliminaries}\label{sec2}
In this section, we recall some facts about  the first and second Lucas sequences. For a prime number $p$ and a non-zero integer $n$, the $p$-adic valuation of $n$ denoted by $\nu_p(n)$ is the exponent of $p$ in the factorization of $n$.
Sanna \cite{Sanna} derived the following formulae for the $p$-adic valuation of the first  Lucas sequence $(U_n)_{n\geq0}$.
\begin{lemma}{\rm(\cite{Sanna}, Corollary 1.6)}\label{lem1}
Let $(U_n)_{n\geq0}$  be the first  Lucas sequence as in \eqref{eq7.15.1} and let $p\geq 3$ be a prime number such that $p\nmid b$. Then for any positive integer $n$, we have
\[
\nu_p(U_n)=\left\{
\begin{aligned}
&\nu _p(n)+\nu _p(U_p)-1, &if~&\, p\mid  \Delta,\,p\mid n, \\
&0, &if~&\, p\mid \Delta,\,p\nmid n, \\
&\nu _p(n)+\nu _p(U_{\tau (p)}),&if~&\, p\nmid \Delta,\, \tau (p)\mid n, \\
&0, &if~&\, p\nmid \Delta,\, \tau (p)\nmid n,\\
\end{aligned}
\right.
\]
where $\Delta=a^2+4b$.
\end{lemma}

\begin{lemma}{\rm (\cite{Sanna}, Theorem 1.5 for $p=2$) }\label{lem2}
Let $(U_n)_{n\geq0}$ be the first  Lucas sequence as in \eqref{eq7.15.1}. If $2\nmid b$, then for any positive integer $n$, we have
\[
\nu_2(U_n)=\left\{
\begin{aligned}
&\nu _2(n)+\nu _2(a)-1, &if~&\,   2\mid a,\, 2\mid n,\\
&0,&if~&\,  2\mid a,\, 2\nmid n, \\
&\nu_2(n)+\nu_2(U_6)-1, &if~&\,  2\nmid a,\,3\mid n,\, 2\mid n,\\
&\nu_2(U_3), &if~&\, 2\nmid a,\, 3\mid n,\, 2\nmid n, \\
&0, &if~&\, 2\nmid a,\,3\nmid n.
\end{aligned}
\right.
\]
\end{lemma}

 Onphaeng et al.\cite{Onphaeng} derived the following formulae for the $p$-adic valuation of the second  Lucas sequence $(V_n)_{n\geq0}$.
\begin{lemma}{\rm (\cite{Onphaeng}) }\label{lemma:3.312}
Let $(V_n)_{n\geq0}$ be the second  Lucas sequence as in \eqref{eq7.15.2} and let $p\geq3$ be a prime number such that $p\nmid b$. Then we have
\[
\nu_p(V_n)=\left\{
\begin{aligned}
&\nu_p(n)+\nu_p(U_{\tau(p)}), &&if\,p\nmid \Delta,\,  \tau (p)\nmid n,\, \tau (p)\mid 2n,\\
&0, &&\text{otherwise}.\\
\end{aligned}
\right.
\]
\end{lemma}
\begin{lemma}{\rm (\cite{Onphaeng}) }\label{lemma:3.311}
Let $(V_n)_{n\geq0}$  be the second  Lucas sequence as in \eqref{eq7.15.2}. If $2\nmid b$, then for any positive integer $n$, we have
\[
\nu_2(V_n)=\left\{
\begin{aligned}
&1, &if&\,        2\mid a,\, 2\mid n,\\
&\nu_2(a),&if&\,  2\mid a,\, 2\nmid n, \\
&1, &if&\,  2\nmid a,\,3\mid n,\, 2\mid n,\\
&\nu_2(a^2+3b), &if&\, 2\nmid a,\, 3\mid n,\, 2\nmid n, \\
&0, &if&\, 2\nmid a,\,3\nmid n.
\end{aligned}
\right.
\]
\end{lemma}

McDaniel \cite{McDaniel} presented the following results regarding the greatest common divisor of the first and second Lucas sequences.
\begin{lemma}{\rm (\cite{McDaniel}, Main Theorem)}\label{lem 3}
Let $(U_n)_{n\geq0}$ and $(V_n)_{n\geq0}$ be the first and second Lucas sequences as in \eqref{eq7.15.1} and \eqref{eq7.15.2}, respectively.  Let $m,\,n \geq 3$ be integers and let $d=\gcd(m,\,n)$. Then the following statements hold.
\begin{enumerate}
\item[\rm(i)]\,$\gcd(U_m,\,U_n)=U_d$;

\item[\rm(ii)]\, $\gcd(V_m,\,V_n)=\left\{
\begin{aligned}
                   &V_d,\  &&if\  \nu_2(m)=\nu_2(n),\\
                   &1 \ \text{or}\  2,\  &&\text{otherwise};
 \end{aligned}
                   \right.
                   $

\item[\rm(iii)]\,  $\gcd(U_m,\,V_n)=\left\{
\begin{aligned}
                   &V_d,\  &&if\  \nu_2(m)> \nu_2(n),\\
                   &1\  or\  2,\  && \text{otherwise}.
\end{aligned}
                   \right.
                   $
\end{enumerate}
\end{lemma}
\begin{lemma}{\rm (\cite{McDaniel})}\label{lem 7.1}
Let $(U_n)_{n\geq0}$ and $(V_n)_{n\geq0}$ be the first and second Lucas sequences as in \eqref{eq7.15.1} and \eqref{eq7.15.2}, respectively.  Let $m,\,n \geq 3$ be integers and let $d=\gcd(m,\,n)$. Then the following statements hold.
\begin{enumerate}
\item[\rm(i)]\,  $2\mid \gcd(U_m,\,V_n)$ if and only if  $2 \mid a$,  $2 \mid m$ or $2 \nmid a$, $2 \nmid b$,\  $3 \mid d$.
 \item[\rm(ii)]\, $2\mid \gcd(V_m,\,V_n)$ if and only if  $2 \mid a$ or $2 \nmid a$, $2 \nmid b$,\  $3 \mid d$.
\end{enumerate}
\end{lemma}

The next lemma summarizes some basic divisibility properties of
the first and second Lucas sequences.

\begin{lemma}{\rm (\cite{Renault, Sanna})}\label{lem 7.24.1}
Let $(U_n)_{n\geq0}$ and $(V_n)_{n\geq0}$ be the first and second Lucas sequences as in \eqref{eq7.15.1} and \eqref{eq7.15.2}, respectively.
Let $p$ be a prime number.
 Then the following statements hold.
\begin{enumerate}
\item[\rm(i)]\,  If $n\geq1$, then $\gcd(U_n,\,b)=1$ and $\gcd(V_n,\,b)=1$.
\item[\rm(ii)]\, For each integer $m\geq1$,  $ m\mid U_n$ for some positive integer $n$ if and only if $\gcd(m,\,b)=1$ and $\tau(m)\mid n$.
\item[\rm(iii)]\, If $p\nmid b$, then $\tau(p)=p$
if and only if $p\mid \tau(p)$ if and only if $p\mid \Delta$, where $\Delta=a^2+4b$.
\end{enumerate}
\end{lemma}
Next, we give the monotonicity of the first and second Lucas sequences  as follows.
\begin{lemma}\label{lem7.6.1}
Let $(U_n)_{n\geq0}$ and $(V_n)_{n\geq0}$ be the first and second Lucas sequences as in \eqref{eq7.15.1} and \eqref{eq7.15.2}, respectively.  If $\Delta=a^2+4b>0$, then for any positive integer $n$, we have $|U_{n+2}|>|U_{n+1}|$ and $|V_{n+1}|>|V_{n}|$.
\end{lemma}

\begin{prof}
If $a\geq 1$ and $b>0$, then it follows from \eqref{eq7.15.1} and \eqref{eq7.15.2} that
\[
U_{n+2}-U_{n+1}=(aU_{n+1}+bU_{n})-U_{n+1}=(a-1)U_{n+1}+bU_{n}>0
\]
for any positive integer $n$. Similarly,
 we have
\[
V_{n+1}-V_{n}=(aV_{n}+bV_{n-1})-V_{n}=(a-1)V_{n}+bV_{n-1}>0
\]
for any positive integer $n$.
If $a\geq 1$ and $b<0$, then $a=\alpha +\beta \geq 1$ and $-b=\alpha \beta>0$. Since $\Delta>0$, we have  $\alpha >\beta >0$.
Since
\[
\alpha-\beta=\frac{a+\sqrt{\Delta}}{2}-\frac{a-\sqrt{\Delta}}{2}=\sqrt{\Delta}\geq1,
\]
we have $\alpha>1$.
It follows that $\alpha ^n>\beta ^n$ for any positive integer $n$.
By Binet's formulae \eqref{eq1} and \eqref{eq2}, we get
\begin{align}
U_{n+2}-U_{n+1}  &=\frac{\alpha^{n+2}-\beta^{n+2}}{\alpha-\beta}-\frac{\alpha^{n+1}-\beta^{n+1}}{\alpha-\beta}  \notag\\
                   &=\frac {\alpha^{n+1}(\alpha-1)+\beta^{n+1}(1-\beta)}{\alpha-\beta} \notag\\
                   &>\frac {\beta^{n+1}(\alpha-1)+\beta^{n+1}(1-\beta)}{\alpha-\beta} \notag\\
                   &=\frac {\beta^{n+1}(\alpha-\beta)}{\alpha-\beta}>0  \notag
\notag \end{align}
for any positive integer $n$.
Similarly, we have
\begin{align}
V_{n+1}-V_n  &=(\alpha^{n+1}+\beta^{n+1})-(\alpha^{n}+\beta^{n})  \notag\\
                   &=\alpha^n(\alpha-1)+\beta^n(\beta-1)  \notag\\
                   &>\beta^n(\alpha-1)+\beta^n(\beta-1)\notag\\
                  &=\beta^n(\alpha+\beta-2).\notag
\notag \end{align}
Since $\alpha >1$ and $\beta>0$, we have $\alpha+\beta>1$.
If  $\alpha+\beta>2$, then  $V_{n+1}>V_n$.
If $\alpha+\beta=a=2$, then  there does not exist an  integer $b<0$ such that $\Delta=a^2+4b>0$.
Now, for $a\geq1$, we conclude that $U_{n+2}>U_{n+1}$  and $V_{n+1}>V_{n}$ for any positive integer $n$.
Note that by Binet's formulae \eqref{eq1} and \eqref{eq2},
\[
U_n(-a,b)=(-1)^{n-1}U_n(a,b)
\]
and
\[
V_n(-a,b)=(-1)^{n}V_n(a,b).
\]
Therefore, for $a< 0$, we conclude that $|U_{n+2}|>|U_{n+1}|$ and $|V_{n+1}|>|V_{n}|$ for any positive integer $n$.
\end{prof}

The following lemma is very useful in the proof of our theorems.
\begin{lemma}\label{lem5}
Let $(U_n)_{n\geq0}$ and $(V_n)_{n\geq0}$ be the first and second Lucas sequences as in \eqref{eq7.15.1} and \eqref{eq7.15.2}, respectively.
Let $m,\,n \geq 3$ be integers.  If $\Delta=a^2+4b>0$, then we have the following statements.
\begin{enumerate}
\item[\rm(i)]\,   $U_n\mid U_m$ if and only if $n\mid m$.

\item[\rm(ii)]\, $V_n\mid U_m$ if and only if $n\mid m$ and $m/n$ is even.

\item[\rm(iii)]\,  If $U_n\mid m$, then $n\mid \tau (m)$.

\item[\rm(iv)]\,  If $V_n\mid m$, then $2n\mid \tau (m)$.
\end{enumerate}
\end{lemma}

\begin{prof}
\begin{enumerate}

\item[\rm(i)]\,
 We now prove necessity. It follows from Lemma \ref{lem 3} (i) that $\gcd (U_m,\,U_n)=U_{d}$, where $d=\gcd(m,\,n)$. Since $U_n\mid U_m$, we have $U_{d}=U_n$. By Lemma \ref{lem7.6.1}, we have $d=n$. So $n \mid m$.
We next prove sufficiency. Since $n \mid m$,  we obtain $d=n$. Therefore, $\gcd (U_m,\,U_n)=U_{d}=U_n$. Hence,  $U_n\mid U_m$.

\item[\rm(ii)]\, We now prove sufficiency. Since $n\mid m$ and $m/n$ is even, we have $\nu_2(m)> \nu_2(n)$. It follows from Lemma \ref{lem 3} (iii) that $\gcd(U_m,\,V_n)=V_{d}=V_n$, where $d=\gcd(m,\,n)$. So $V_n\mid U_m$. We next prove necessity.
If $\nu_2(m)\leq\nu_2(n)$, then $\gcd(U_m,\,V_n)=V_n=1$ or $2$. By Lemma \ref{lem7.6.1}, we know that $V_n \neq 1$ or $2$ for each integer $n\geq3$.
If $\nu_2(m)>\nu_2(n)$, then $\gcd(U_m,\,V_n)=V_{d}=V_n$. From Lemma \ref{lem7.6.1}, it yields $d=n$.
Thus, $n\mid m$ and $m/n$ is even.

\item[\rm(iii)]\,
Since $U_n\mid m$ and $m\mid U_{\tau (m)}$, we know that $U_n\mid U_{\tau (m)}$.  By Lemma \ref{lem5} (i), we get $n\mid \tau (m)$.

\item[\rm(iv)]\, Since $V_n\mid m$ and $m\mid U_{\tau (m)}$, we see that $V_n\mid U_{\tau (m)}$. By Lemma \ref{lem5} (ii), we obtain that $\tau (m)/n$ is even. In particular, $2n\mid \tau (m)$.
\end{enumerate}
\end{prof}
\begin{remark}
\rm
 Lemma \ref{lem5} (i) and (ii) fail to hold for $\Delta=a^2+4b<0$. Therefore, we impose the condition that $\Delta=a^2+4b>0$.
 Some counterexamples are presented below.
 For $a=-3$ and $b=-5$, we have $U_4=3\mid U_6=72$ but $4 \nmid 6$.
 For $a=1$ and $b=-2$, we have $U_8=-3\mid U_{12}=45$ but $8 \nmid 12$.
 For $a=4$ and $b=-5$, we have $V_3=4\mid U_{4}=24$ but $3 \nmid 4$.
 For $a=2$ and $b=-3$, we have $V_5=2\mid U_{6}=-10$ but $5 \nmid 6$.
\end{remark}

By the recursive formula for $(U_n)_{n\geq0}$ as in \eqref{eq7.15.1}, we get the following result.

\begin{lemma}\label{lem 7.13.3}
Let $(U_n)_{n\geq0}$  be the first  Lucas sequence as in \eqref{eq7.15.1}. If  $a$  is   odd, then we have the following statements.
\begin{enumerate}
\item[\rm(i)]\,  If $b\equiv 1 \ (\rm mod\ 4)$, then $\nu_2(U_3)=1$ and $\nu_2(U_6)=\nu_2(a^2+3b)+1$.

\item[\rm(ii)]\,  If $b\equiv 3\  (\rm mod\ 4)$, then $\nu_2(U_3)\geq 2$ and $\nu_2(U_6)=\nu_2(U_3)+1$.
\end{enumerate}
\end{lemma}

\begin{prof}
\begin{enumerate}

\item[\rm(i)]\,
Since $U_3=a^2+b\equiv 2 \ (\rm mod\ 4)$, we have $\nu_2(U_3)=1$.
By the recursive formula for $(U_n)_{n\geq0}$, we know that $U_6=aU_3(a^2+3b)$.
Therefore, $\nu_2(U_6)=\nu_2(a^2+3b)+1$.

\item[\rm(ii)]\, Since $a^2+3b\equiv 2\ (\rm mod\ 4)$, we have  $\nu_2(a^2+3b)=1$.
Since $U_3=a^2+b\equiv 0\ (\rm mod\ 4)$,  we get $\nu_2(U_3)\geq 2$.
Therefore, $\nu_2(U_6)=\nu_2(U_3)+1$.

\end{enumerate}
\end{prof}

\section{The Proof of Theorems}\label{sec3}
 In this section, we will give the proofs of  Theorem \ref{Main theorem} and Theorem \ref{Th1.2}.

\textbf{The proof of Theorem 1.1 (i)}
\begin{prof}
Since $U_m\mid U_m V_n$, it follows from  Lemma \ref{lem5} (iii) that $m\mid \tau (U_m V_n)$.
Since $V_n\mid U_m V_n$, it follows from Lemma \ref{lem5} (iv) that  $2n\mid \tau (U_m V_n)$.
Hence, we obtain
\begin{equation}\label{eq3}
[m,2n]\mid \tau (U_m V_n).
\end{equation}
By Lemma \ref{lem 3} (iii), we know that $\gcd(U_m,\,V_n)=1$, $2$ or $V_d$, where $d=\gcd(m,\, n)$.
The proof will be presented in two cases based on $\nu_2(m) \leq \nu_2(n)$ and  $\nu_2(m) > \nu_2(n)$.

 \textbf {Case 1.}
If $\nu_2(m) \leq \nu_2(n)$, then  $\gcd(U_m,\,V_n)=1$ or $2$. Thus $[m,2n]=2[m,n]$.
From \eqref{eq3}, it follows that
\begin{equation}\label{eq7.8.1}
 2[m,n]\mid \tau(U_mV_n).
 \end{equation}
Lemma \ref{lem5} (i) and (ii) yield that $U_m\mid  U_{2[m,n]}$ and $V_n\mid  U_{2[m,n]}$.
If $\gcd(U_m,\,V_n)=1$, then $U_m V_n \mid  U_{2[m,n]}$. Therefore, by Lemma \ref{lem 7.24.1} (ii), we obtain
$\tau (U_m V_n)\mid 2[m,n]$.
Together with \eqref{eq7.8.1}, we have  $\tau (U_m V_n)=2[m,n]$ as claimed.
Next we can assume that $\gcd(U_m,\,V_n)=2$.
It suffices to show that $\nu_2(U_mV_n)\leq \nu_2(U_{2[m,n]})$.
By Lemma \ref{lem 7.1} (i), we divide the calculation into three subcases.

\textbf {Subcase 1.1.}  $2\mid a$, $2\mid m$.
Since $\nu_2(m)\leq \nu_2(n)$, we have $2\mid n$.
From Lemmas \ref{lem2} and  \ref{lemma:3.311}, it follows that
\begin{align}
\nu _2\left(U_{2[m,n]}\right) &=1+\nu _2([m,n])+\nu _2(a)-1  \notag\\
                   &\geq \nu _2(m)+\nu _2(a)-1+1  \notag\\
                   &=\nu_2(U_m)+\nu_2(V_n)  \notag\\
                   &=\nu_2(U_mV_n).  \notag
\notag \end{align}

\textbf {Subcase 1.2.} $2\nmid b$, $2\nmid a$, $3\mid d$, $2\mid m$.
Since $\nu_2(m)\leq \nu_2(n)$, we have $2\mid n$.
Thus, we get
\begin{align}
\nu _2\left(U_{2[m,n]}\right) &=1+\nu _2([m,n])+\nu _2(U_6)-1  \notag\\
                   &\geq \nu _2(m)+\nu _2(U_6)-1+1  \notag\\
                   &=\nu_2(U_m)+\nu_2(V_n)  \notag\\
                   &=\nu_2(U_mV_n).  \notag
\notag \end{align}

\textbf {Subcase 1.3.} $2\nmid b$, $2\nmid a$, $3\mid d$, $2\nmid m$.
By Lemma \ref{lemma:3.311},  $\nu_2(V_n)$ is equal to either $1$ or $\nu_2(a^2+3b)$.
Since $\nu_2(a^2+3b)\geq 1$, we have $\nu_2(V_n)\leq \nu_2(a^2+3b)$. Using the identity $U_6=a(a^2+3b)U_3$, we get
\begin{align}
\nu _2\left(U_{2[m,n]}\right) &=1+\nu _2([m,n])+\nu _2(U_6)-1  \notag\\
                   &= \nu _2([m,n])+\nu _2(U_3)+\nu_2(a^2+3b)  \notag\\
                   &\geq \nu_2(U_3)+\nu_2(a^2+3b)  \notag\\
                   &= \nu_2(U_m)+\nu_2(a^2+3b)  \notag\\
                   &\geq \nu_2(U_m)+\nu_2(V_n)  \notag\\
                   &=\nu_2(U_mV_n).  \notag
\notag \end{align}
Therefore, we obtain $\tau (U_m V_n)=2[m,n]$ as desired.

\textbf {Case 2.}
If $\nu_2(m)> \nu_2(n)$, then $\gcd(U_m,\,V_n)=V_d$.
If $\gcd(U_m,\,V_n)=V_d=1$, then $U_m V_n \mid U_{[m, n]}$. So, $\tau(U_m V_n) \mid [m, n]$.
By \eqref{eq3}, we obtain  $[m, n] \mid \tau(U_mV_n)$. Therefore, $\tau(U_m V_n) = [m, n]$.
Suppose $\gcd(U_m,\,V_n)=V_d\geq2$,
and let $p$ be  any prime factor of $V_d$. By Lemma \ref{lem 7.24.1} (i), we know that  $p\nmid b$. Thus, $\nu_p(U_m)>0$ and $\nu_p(V_n)>0$.
 If $p\geq 3$, then it follows from  Lemmas \ref{lem1} and \ref{lemma:3.312} that $p\nmid \Delta$ and $\tau (p)\mid m $.
Thus, $\tau (p)\mid [m,n]$.
It follows that
\begin{align}\label{eq.9.7.1}
 \nu_p\left(U_{[m,n] V_d}\right)&=\nu_p([m,n]V_d)+\nu_p(U_{\tau(p)})\notag\\
                         &=\nu_p([m,n])+\nu_p(U_{\tau(p)})+\nu_p(d)+\nu_p(U_{\tau(p)})\notag\\
                        &=\nu_p(m)+\nu_p(U_{\tau(p)})+\nu_p(n)+\nu_p(U_{\tau(p)})\notag\\
                      &=\nu_p(U_m)+\nu_p(V_n)\notag\\
                       &=\nu_p(U_m V_n).
 \end{align}
 If $p=2$, then $2\mid V_d$. By Lemma \ref{lemma:3.311}, we consider the following two subcases. Since $\nu_2(m)> \nu_2(n)$, we have $2\mid m$ and $2 \mid [m,n]$.

\textbf {Subcase 2.1.} $  2\mid a$.
By Lemma \ref{lemma:3.311}, $\nu_2(V_d)=\nu_2(V_n)$ holds irrespective  of whether $2\mid d$ or $2\nmid d$.
Therefore, we obtain
\begin{align}\label{eq.9.7.2}
\nu_2\left(U_{[m,n] V_d}\right) &=\nu _2([m,n] V_d)+\nu _2(a)-1  \notag\\
                   &=\nu _2(m)+\nu_2(V_n)+\nu _2(a)-1  \notag\\
                   &= \nu_2(U_m)+\nu_2(V_n)  \notag\\
                   &=\nu_2(U_mV_n).
 \end{align}

\textbf {Subcase 2.2.} $2\nmid b,\  2\nmid a,\  3\mid d$.
By Lemma \ref{lemma:3.311}, $\nu_2(V_d)=\nu_2(V_n)$ holds irrespective  of whether $2\mid d$ or $2\nmid d$.
Hence, we get
\begin{align}\label{eq.9.7.3}
\nu_2\left(U_{[m,n] V_d}\right)&=\nu _2([m,n]V_d)+\nu _2(U_6)-1 \notag\\
                    &=\nu _2(m)+\nu _2(V_n)+\nu _2(U_6)-1 \notag\\
                   &= \nu_2(U_m)+\nu_2(V_n)  \notag\\
                   &=\nu_2(U_mV_n).
\end{align}
Lemma \ref{lem5} (i) and (ii) yield that $U_m\mid  U_{[m,n]V_d}$ and $V_n\mid  U_{[m,n]V_d}$.
By \eqref{eq.9.7.1}--\eqref{eq.9.7.3}, we know that $U_mV_n\mid  U_{[m,n]V_d}$. Thus, $\tau(U_m V_n) \mid [m,n]V_d$.
From \eqref{eq3}, we have $[m,n]\mid \tau(U_mV_n)$.
Now, we claim that
\[
\tau(U_m V_n)=[m,n] V_d.
 \]
 Let us assume the opposite, namely that there exists an integer $t$ such that
$\tau(U_m V_n) = [m,n]t $, where $t \mid V_d$ but $t \neq V_d$.
It follows that there exists a prime number $q$ such that $\nu_{q}(V_d)> \nu_{q}(t)$.
So, $q$ is the prime factor of $V_d$.
From the previous discussion, we conclude that
\[
\nu_{q}\left(U_{[m,n]t}\right)<\nu_{q}\left(U_{[m,n] V_d}\right)=\nu _{q}(U_m V_n),
\]
which is a contradiction. The proof is then complete.
\end{prof}

\textbf{The proof of Theorem 1.1 (ii)}
\begin{prof}
We know that $\gcd(U_m,\,U_n)=U_d$, where $d=\gcd(m,\, n)$.
Since $U_m\mid U_m U_n$ and $U_n\mid U_m U_n$,  we obtain that $m\mid \tau (U_m U_n)$ and $n\mid \tau (U_m U_n)$.
So, we have
\begin{equation}\label{eq4}
[m,n] \mid \tau (U_m U_n).
\end{equation}
By Lemma \ref{lem5} (i), we know that $U_m \mid U_{[m,n]}$ and $U_n\mid  U_{[m,n]}$.

If $\gcd(U_m,\,U_n)=1$,  then $U_m U_n\mid U_{[m,n]}$. Therefore, by Lemma \ref{lem 7.24.1} (ii), we obtain
 $\tau (U_m U_n) \mid [m,n]$. Together with \eqref{eq4}, we obtain $\tau (U_m U_n) = [m,n]$. From now on, we can assume that $\gcd(U_m,\,U_n)=U_d\geq 2$. Let us consider any prime factor $p$ of $U_d$. Then $p\nmid b$.
For $p\geq 3$, it follows from Lemma \ref{lem1} that either $p\mid \Delta$, $p\mid d$  or $p\nmid \Delta$, $\tau(p)\mid d$.
If $p\mid \Delta$ and $p\mid d$, then $p\mid [m,n]$. Thus, we get
\begin{align}\label{eq7.9.1}
 \nu_p\left(U_{[m,n] U_d}\right)&=\nu_p([m,n]U_d)+\nu_p(U_p)-1\notag\\
                      &=\nu_p([m,n])+\nu_p(d)+\nu_p(U_p)-1+\nu_p(U_p)-1\notag\\
                        &=\nu_p(m)+\nu_p(n)+\nu_p(U_p)-1+\nu_p(U_p)-1\notag\\
                      &=\nu_p(U_m)+\nu_p(U_n)\notag\\
                       &=\nu_p(U_m U_n).
 \end{align}
 If $p\nmid \Delta$ and $\tau(p)\mid d$, then $\tau(p)\mid [m,n]$. So, we have
\begin{align}\label{eq7.9.2}
 \nu_p\left(U_{[m,n] U_d}\right)&=\nu_p([m,n]U_d)+\nu_p(U_{\tau(p)})\notag\\
                   &=\nu_p([m,n])+\nu_p(d)+2\nu_p(U_{\tau(p)})\notag\\
                        &=\nu_p(m)+\nu_p(n)+2\nu_p(U_{\tau(p)})\notag\\
                       &=\nu_p(U_m U_n).
 \end{align}
If $p=2$, then $2\mid U_d$. By Lemma \ref{lem2}, we consider the following three cases.

\textbf {Case 1.}  $2\mid a$, $2\mid d$. Then $2\mid [m,n]$.
We have
\begin{align}
\nu_2\left(U_{[m,n] U_d}\right) &=\nu _2([m,n]U_d)+\nu _2(a)-1\notag\\
             &=\nu _2([m,n])+\nu_2(d)+\nu _2(a)-1+\nu _2(a)-1  \notag\\
                   &= \nu_2(U_m)+\nu_2(U_n) \notag\\
                   &=\nu_2(U_mU_n). \notag
\notag \end{align}

\textbf {Case 2.} $2\nmid b$, $2\nmid a$, $3\mid d$, $2\mid d$.
Then $2\mid [m,n]$. We get
\begin{align}
\nu_2\left(U_{[m,n] U_d}\right) &=\nu _2([m,n]U_d)+\nu _2(U_6)-1\notag\\
&=\nu _2([m,n])+\nu _2(d)+\nu _2(U_6)-1+\nu _2(U_6)-1\notag\\
                   &= \nu_2(U_m)+\nu_2(U_n)  \notag\\
                   &=\nu_2(U_mU_n). \notag
\notag \end{align}

\textbf {Case 3.} $2\nmid b$, $2\nmid a$, $3\mid  d$, $2\nmid d$.
In this case, we have $\nu_2(U_d)=\nu_2(U_3)$.
We consider three subcases according to the parity of $m$ and $n$.

\textbf {Subcase 3.1.} $2\nmid b$, $2\nmid a$, $3\mid  d$, $2\nmid m$, $2\mid n$. Then $3\mid [m,n]$ and $2\mid [m,n]$. We have
\begin{align}
\nu_2\left(U_{[m,n] U_d}\right)
&=\nu _2([m,n]U_d)+\nu _2(U_6)-1\notag\\
 &=\nu _2([m,n])+\nu _2(U_3)+\nu _2(U_6)-1\notag\\
                   &=\nu _2(n)+\nu _2(U_6)-1+\nu _2(U_3)\notag\\
                   &= \nu_2(U_n)+\nu_2(U_m)  \notag\\
                   &=\nu_2(U_mU_n). \notag
\notag \end{align}

\textbf {Subcase 3.2.} $2\nmid b$, $2\nmid a$, $3\mid  d$, $2\mid m$, $2\nmid n$.  It is  similar to Subcase 3.1 and we obtain
\[
\nu_2\left(U_{[m,n] U_d}\right)=\nu_2(U_mU_n).
\]

\textbf {Subcase 3.3.} $2\nmid b$, $2\nmid a$, $3\mid  d$, $2\nmid m$, $2\nmid n$.
Thus, $\nu_2(U_mU_n)=2\nu_2(U_3)$.
If  $b\equiv 1\ (\rm mod\ 4)$, then it follows from Lemma \ref{lem 7.13.3} (i) that
\begin{align}
\nu_2\left(U_{[m,n] U_d}\right) &=\nu _2([m,n])+\nu _2(U_3)+\nu _2(U_6)-1\notag\\
                  &=0+1+\nu _2(U_6)-1\notag\\
                   &=\nu _2(a^2+3b)+1\notag\\
                   &\geq 2=2\nu_2(U_3)\notag\\
                   &=\nu_2(U_mU_n).\notag
\notag \end{align}
If $b\equiv 3\ (\rm mod\ 4)$, then it follows from Lemma \ref{lem 7.13.3} (ii) that
\begin{align}
\nu_2\left(U_{[m,n]U_d}\right) &=\nu _2([m,n])+\nu _2(U_3)+\nu _2(U_6)-1\notag\\
                         &=0+\nu _2(U_3)+(\nu _2(U_3)+1)-1\notag\\
                       &=2\nu _2(U_3)\notag\\
                   &=\nu_2(U_mU_n).\notag
\notag \end{align}
Lemma \ref{lem5} (i)  yields that $U_m\mid  U_{[m,n]U_d}$ and $U_n\mid  U_{[m,n]U_d}$.
From above discussion, we can conclude that $\nu _p(U_m U_n) \leq \nu_p\left(U_{[m,n] U_d}\right)$ for any prime factor $p$ of $U_d$.
So, $U_m U_n \mid U_{[m,n]U_d}$, yielding
\[
\tau(U_m U_n)\mid [m,n] U_d.
\]
Now, together with \eqref{eq4}, we claim that
\[
\tau(U_m U_n)=[m,n] U_d.
\]
Let us assume the opposite, namely that there exists an integer $t$ such that $\tau(U_m U_n)=[m,n]t $, where $t\mid U_d$ but $t\neq U_d$.
It follows that there exists a prime number $q$ such that $\nu_{q}(U_d)> \nu_{q}(t)$.
So, $q$ is the prime factor of $U_d$.
If there exists an odd prime number $q$ such that $\nu_{q}(U_d)> \nu_{q}(t)$, then it follows from \eqref{eq7.9.1} and \eqref{eq7.9.2} that $\nu_{q}(U_{[m,n]t })<\nu _{q}(U_m U_n)$, which is a contradiction.
Otherwise, $\nu_{2}(U_d)> \nu_{2}(t)$ and $\nu_{q}(U_d)= \nu_{q}(t)$ for every odd prime number $q$.
We only  need to consider the case when $2\nmid b$, $2\nmid a$, $3\mid d$, $2\nmid m$, $2\nmid n$.
If $t=U_d/2$ is odd,  then
\[
 \nu_2(U_{[m,n] t })=\nu_2(U_3)<2\nu_2(U_3)=\nu _2(U_m U_n),
\]
which is a contradiction.
If $t=U_d/2$ is even, then it suffices to show $U_m U_n \nmid U_{[m,n] U_d/2}$.
Note that $\nu_2(U_d)=\nu_2(U_3)\geq 2$. By Lemma \ref{lem 7.13.3} (ii), we have $b\equiv 3\ (\rm mod\ 4)$.
Thus
\begin{align}
\nu_2\left(U_{[m,n] U_d/2}\right) &=\nu _2([m,n])+\nu_2(U_d/2)+\nu_2(U_6)-1\notag\\
                           &=0+\nu_2(U_3)-1+\nu_2(U_6)-1\notag\\
                           &=\nu_2(U_3)-1+(\nu_2(U_3)+1)-1\notag\\
                           &=2\nu_2(U_3)-1\notag\\
                   &<\nu_2(U_mU_n),\notag
\notag \end{align}
and the claim follows.
\end{prof}

\textbf{The proof of Theorem 1.1 (iii)}
\begin{prof}
Since $V_m\mid V_m V_n$ and $V_n\mid V_m V_n$, we obtain that $2m\mid \tau (V_m V_n)$ and $2n\mid \tau (V_m V_n)$. So, we have
\begin{equation}\label{eq5}
2[m,n]\mid \tau (V_m V_n).
\end{equation}
By Lemma \ref{lem 3} (ii), we know that
$\gcd(V_m,\,V_n)=1$, $2$ or $V_d$, where $d=\gcd(m,\, n)$.
The proof will be presented in two cases based on $\nu_2(m) \neq \nu_2(n)$ and  $\nu_2(m) = \nu_2(n)$.

\textbf {Case 1.}
If $\nu_2(m)\neq \nu_2(n)$, then $\gcd(V_m,\,V_n)=1$ or $2$.
By Lemma \ref{lem5} (ii), we know that $V_m \mid  U_{2[m,n]}$ and $V_n\mid  U_{2[m,n]}$.
If $\gcd(V_m,\,V_n)=1$, then  $V_m V_n\mid  U_{2[m,n]}$.
Therefore, we obtain $\tau (V_m V_n)\mid 2[m,n]$.
Together with \eqref{eq5}, it follows that  $\tau (V_m V_n)=2[m,n]$.
Next we can assume that  $\gcd(V_m,\,V_n)=2$. Then we have $\nu_2(m)\neq \nu_2(n)$.  Hence, at least one of the numbers $m$ and $n$ is even. It is enough to show that $\nu_2(V_m V_n)\leq \nu_2(U_{2[m,n]})$.
By Lemma \ref{lem 7.1} (ii), we only need to consider the following two subcases.

\textbf {Subcase 1.1.}  $2\mid a$.
By Lemma \ref{lemma:3.311}, we get $\nu_2(V_mV_n)\leq \nu_2(a)+1$.
Therefore, we have
\begin{align}
\nu _2\left(U_{2[m,n]}\right) &=1+\nu _2([m,n])+\nu _2(a)-1  \notag\\
                   &=\max\{\nu_2(m),\nu_2(n)\}+\nu_2(a)  \notag\\
                   &\geq \nu_2(a)+1  \notag\\
                   &\geq \nu_2(V_mV_n).  \notag
\notag \end{align}

\textbf {Subcase 1.2.} $2\nmid b$, $2\nmid a$, $3\mid d$.
By Lemma \ref{lemma:3.311}, we get $\nu_2(V_mV_n)\leq \nu_2(a^2+3b)+1$.
Therefore, we have
\begin{align}
\nu _2\left(U_{2[m,n]}\right) &=1+\nu _2([m,n])+\nu _2(U_6)-1  \notag\\
                   &=\max\{\nu_2(m),\nu_2(n)\}+\nu_2(a^2+3b)+\nu_2(U_3)  \notag\\
                   &\geq \nu_2(a^2+3b)+1 \notag\\
                   &\geq \nu_2(V_mV_n).  \notag
\notag \end{align}
It follows that $V_mV_n\mid U_{2[m,n]}$, yielding $\tau(V_mV_n)\mid 2[m,n]$.
Together with \eqref{eq5}, we obtain $\tau (V_m V_n)=2[m,n]$.

\textbf {Case 2.}
If $\nu_2(m)=\nu_2(n)$, then $\gcd(V_m,\,V_n)=V_d$.
If $\gcd(V_m, V_n) = V_d=1$, then it follows from Lemma \ref{lem7.6.1} that  $d=1$ and $a=1$.
Since $V_m \mid U_{2[m,n]}$ and $V_n \mid U_{2[m,n]}$,
it follows that $V_m V_n \mid U_{2[m,n]}$.
Therefore, we obtain $\tau(V_m V_n) \mid 2[m,n]$.
By \eqref{eq5}, we conclude that $\tau(V_m V_n) = 2[m,n]$.
Suppose $\gcd(V_m,\,V_n)=V_d\geq 2$, and
let $p$ be  any prime factor of $V_d$. Then $p\nmid b$,
$\nu_p(V_m)>0$ and $\nu_p(V_n)>0$.
If $p\geq 3$,  then it follows from Lemma \ref{lemma:3.312} that
$p\nmid \Delta$,  $\tau (p)\nmid m $, $\tau (p)\mid 2m $, $\tau (p)\nmid n$ and $\tau (p)\mid  2n $.
So, $\tau (p)\nmid [m,n]$ and $\tau (p)\mid 2[m,n]$.
Thus, we have
\begin{align}
 \nu_p\left(U_{2[m,n]V_d}\right)&=\nu_p(2[m,n]V_d)+\nu_p(U_{\tau(p)})\notag\\
                        &=\nu_p(2[m,n])+\nu_p(d)+\nu_p(U_{\tau(p)})+\nu_p(U_{\tau(p)})\notag\\
                        &=\nu_p(m)+\nu_p(n)+2\nu_p(U_{\tau(p)})\notag\\
                      &=\nu_p(V_m)+\nu_p(V_n)\notag\\
                       &=\nu_p(V_m V_n).\notag
 \notag \end{align}
 It is obvious that $V_m \mid U_{2[m,n]V_d}$ and $V_n\mid U_{2[m,n]V_d}$.
If $V_d$ is odd, then $2\mid b$ or  $2\nmid b,\  2\nmid a, \ 3\nmid d$.
We have $V_m V_n\mid U_{2[m,n]V_d}$.
Assume that $\tau(V_m V_n)=2[m,n]t$, where $t\mid V_d$ but $t\neq V_d$.
It follows that there exists a prime number $q\geq3$ such that $\nu_{q}(V_d)> \nu_{q}(t)$.
So, $q$ is the prime factor of $V_d$.
Therefore, we have
\begin{align}
 \nu_{q}\left(U_{2[m,n] t}\right)
 &=\nu_{q}(2[m,n] t)+\nu_{q}(U_{\tau(q)})\notag\\
  &<\nu_{q}(2[m,n]V_d)+\nu_{q}(U_{\tau(q)})\notag\\
  &=\nu_{q}(V_m V_n).\notag
 \notag \end{align}
Together with the previous discussion, it follows that $\tau(V_m V_n)=2[m,n] V_d$.
 If $V_d$ is even, then by Lemma \ref{lemma:3.311}, we consider the following four subcases.

\textbf {Subcase 2.1.}  $2\mid a$, $2\mid d$.
We have
\begin{align}
\nu_2\left(U_{2[m,n] V_d/2}\right) &=\nu_2\left(U_{[m,n] V_d}\right)\notag\\
                        &=\nu_2([m,n])+1+\nu_2(a)-1\notag\\
                        &=\nu_2(m)+\nu_2(a)\notag\\
                      &\geq 1+1=\nu_2(V_m V_n).\notag
 \notag \end{align}
Note that $\nu_2(V_d)=1$. So, $V_d/4$ is not an integer.
Then   $2[m,n]\nmid 2[m,n] V_d/4$.
Together with \eqref{eq5}, we obtain $\tau(V_m V_n)=[m,n]V_d$.

\textbf {Subcase 2.2.} $2\mid a$, $2\nmid d$.
Since $\nu_2(m)=\nu_2(n)$, we obtain that  $2\nmid m$ and $2\nmid n$. We obtain
\begin{align}
 \nu_2\left(U_{2[m,n] V_d}\right)&=\nu_2(2[m,n] V_d)+\nu_2(a)-1\notag\\
 &=1+0+\nu_2(a)+\nu_2(a)-1\notag\\
                      &=2\nu_2(a) =\nu_2(V_m V_n).\notag
 \notag \end{align}
Now, it is enough to show that $\nu_2(V_m V_n)>\nu_2\left(U_{[m,n] V_d}\right)$. Here
\begin{align}
 \nu_2\left(U_{[m,n] V_d}\right)&=\nu_2([m,n] V_d)+\nu_2(a)-1\notag\\
                        &=2\nu_2(a)-1\notag\\
                      &<\nu_2(V_m V_n).\notag
 \notag \end{align}
Therefore, $\tau(V_m V_n)=2[m,n]V_d$.

\textbf {Subcase 2.3.} $2\nmid b$, $2\nmid a$, $3\mid d$, $2\mid d$.  We get
\begin{align}
 \nu_2\left(U_{[m,n] V_d}\right)&=\nu_2([m,n])+\nu_2(V_d)+\nu_2(U_6)-1\notag\\
                        &=\nu_2(m)+1+\nu_2(U_6)-1\notag\\
                        &=\nu_2(m)+\nu_2(a^2+3b)+\nu_2(U_3)\notag\\
                      &\geq 1+1=\nu_2(V_m V_n).\notag
 \notag \end{align}
 Note that $\nu_2(V_d)=1$. So, $V_d/4$ is not an integer.
Then   $2[m,n]\nmid 2[m,n] V_d/4$.
Together with  \eqref{eq5}, we obtain $\tau(V_m V_n)=[m,n] V_d$.

\textbf {Subcase 2.4.} $2\nmid b$, $2\nmid a$, $3\mid  d$, $2\nmid d$. Then $2\nmid m$ and $2\nmid n$.
If $b\equiv 1\ (\rm mod\ 4)$, then it follows from   Lemma \ref{lem 7.13.3} (i) that
\begin{align}
 \nu_2\left(U_{[m,n] V_d}\right)&=\nu_2([m,n])+\nu_2(V_d)+\nu_2(U_6)-1\notag\\
                        &=0+\nu_2(a^2+3b)+(\nu_2(a^2+3b)+1)-1\notag\\
                        &=2\nu_2(a^2+3b)\notag\\
                      &=\nu_2(V_m V_n).\notag
 \notag \end{align}
In order to draw a conclusion, it suffices to show that $\nu_2\left(U_{2[m,n] V_d/4}\right)<\nu_2(V_m V_n)$.
Note that $\nu_2(V_d)=\nu_2(a^2+3b)\geq2$. Then $V_d/2$ is even.
Hence, we have
\begin{align}
 \nu_2\left(U_{2[m,n]V_d/4}\right)&=\nu_2\left(U_{[m,n]V_d/2}\right)\notag\\
 &=\nu_2([m,n])+\nu_2(V_d/2)+\nu_2(U_6)-1\notag\\
                        &=\nu_2(a^2+3b)-1+\nu_2(U_6)-1\notag\\
                        &=2\nu_2(a^2+3b)-1\notag\\
                      &<\nu_2(V_m V_n).\notag
 \notag \end{align}
 Therefore, $\tau(V_m V_n)=[m,n]V_d$.
 If $b\equiv 3\ (\rm mod\ 4)$, then it follows from Lemma \ref{lem 7.13.3} (ii) that
\begin{align}
 \nu_2\left(U_{[m,n]V_d}\right)&=\nu_2([m,n])+\nu_2(V_d)+\nu_2(U_6)-1\notag\\
                        &=\nu_2(a^2+3b)+(\nu_2(U_3)+1)-1\notag\\
                        &=1+\nu_2(U_3)\notag\\
                        &\geq 2=2\nu_2(a^2+3b)\notag\\
                        &=\nu_2(V_m V_n).\notag
 \notag \end{align}
 Note that $\nu_2(V_d)=\nu_2(a^2+3b)=1$. So, $V_d/4$ is not an integer.
Then   $2[m,n]\nmid 2[m,n] V_d/4$.
Together with \eqref{eq5}, we obtain $\tau(V_m V_n)=[m,n]V_d$.
\end{prof}

\textbf{The proof of Theorem 1.2}
\begin{prof}
We consider the following four cases to complete the proof.

\textbf {Case 1.} $p\nmid n$, $p\geq3$ and $2\nmid n$.
For $s\in \{0,\,p,\,2p\}$, it follows from  Lemma \ref{lem5} (iii) that
\[
(n+s)\mid \tau(U_nU_{n+p}U_{n+2p}).
\]
 Since $\gcd(n,n+p)=1$, $\gcd(n,n+2p)=1$ and $\gcd(n+p,n+2p)=1$,
 we have
\begin{equation}\label{eq7.23.1}
 n(n+p)(n+2p)\mid \tau(U_nU_{n+p}U_{n+2p}).
\end{equation}
 For $s\in \{0,\,p,\,2p\}$, it is obvious that
 \[
 U_{n+s}\mid U_{n(n+p)(n+2p)}.
 \]
 By Lemma \ref{lem 3} (i), we get that
 \[
\gcd(U_n,\, U_{n+p})=1,\  \gcd(U_n,\, U_{n+2p})=1,\
 \gcd(U_{n+p},\, U_{n+2p})=1.
 \]
 This gives
 \[
 U_nU_{n+p}U_{n+2p}\mid U_{n(n+p)(n+2p)}.
 \]
 Thus, we obtain
\begin{equation}\label{eq7.23.2}
 \tau(U_nU_{n+p}U_{n+2p}) \mid n(n+p)(n+2p).
 \end{equation}
Together with \eqref{eq7.23.1} and \eqref{eq7.23.2}, we obtain
\[
\tau(U_nU_{n+p}U_{n+2p})=n(n+p)(n+2p).
\]

\textbf {Case 2.} $p\nmid n$, $p\geq3$ and $2\mid n$.
For $s\in \{0,\,p,\,2p\}$, it follows from  Lemma \ref{lem5} (iii) that
\[
(n+s)\mid \tau(U_nU_{n+p}U_{n+2p}).
\]
It is obvious that
\[
\gcd(n,n+p)=1,\  \gcd(n,n+2p)=2,\  \gcd(n+p,n+2p)=1.
\]
If $\nu_2(n)=1$, then $\gcd(\frac{n}{2},n+2p)=1$.
If $\nu_2(n)\geq 2$, then $\gcd(n,\frac{n+2p}{2})=1$.
Therefore, we get
\begin{equation}\label{eq7.23.3}
\left. \frac{n(n+p)(n+2p)}{2} \ \middle| \ \tau(U_nU_{n+p}U_{n+2p})\right..
\end{equation}
For $s\in \{0,\,p,\,2p\}$, it is obvious that
 \[
 \left. U_{n+s}\ \middle| \ U_{\frac{n(n+p)(n+2p)}{2}}\right..
 \]
Since $\gcd(U_n,\ U_{n+p})=1$, it follows  that $\left. U_nU_{n+p} \ \middle| \ U_{\frac{n(n+p)(n+2p)}{2}}\right.$ and
$\left. U_{n+2p} \ \middle| \ U_{\frac{n(n+p)(n+2p)}{2}}.\right.$
By Lemma \ref{lem 3} (i), we get
\begin{align}
\gcd(U_nU_{n+p},\,U_{n+2p}) &=\gcd(U_n,\,U_{n+2p})\gcd(U_{n+p},\,U_{n+2p}) \notag\\
 &=U_{\gcd(n,\,n+2p)}U_{\gcd(n+p,\,n+2p)}\notag\\
 &=U_2U_1=a. \notag
 \end{align}
If $a=1$, then
\[
\left. U_nU_{n+p}U_{n+2p} \ \middle| \ U_{\frac{n(n+p)(n+2p)}{2}}.\right.
\]
This gives
\begin{equation}\label{eq7.23.4}
\left. \tau(U_nU_{n+p}U_{n+2p}) \ \middle| \ \frac{n(n+p)(n+2p)}{2}.\right.
\end{equation}
Together with \eqref{eq7.23.3} and \eqref{eq7.23.4}, when $a=1$, we obtain
\[
\tau(U_nU_{n+p}U_{n+2p})=\frac{n(n+p)(n+2p)}{2}.
\]
When $a\geq2$,
let us consider any prime factor $q$ of $a$.
It is obvious that
\begin{equation}\label{eq5.13.1}
\left. U_nU_{n+p} \ \middle| \ U_{\frac{n(n+p)(n+2p)}{2}\cdot\frac{a}{\gcd(a,\,n+p)}}\right.
\end{equation}
and
\begin{equation}\label{eq5.13.2}
\left. U_{n+2p} \ \middle| \ U_{\frac{n(n+p)(n+2p)}{2}\cdot\frac{a}{\gcd(a,\,n+p)}}.\right.
\end{equation}
For $q=2$, by Lemma \ref{lem2}, we have
\begin{align}\label{eq8.15.11}
&\nu_2\left(U_{\frac{n(n+p)(n+2p)}{2}\cdot\frac{a}{\gcd(a,\,n+p)}}\right)\notag\\
=&\nu _2(n(n+p)(n+2p))-\nu _2(2)+\nu _2(a)-\nu_2(\gcd(a,\,n+p))+\nu _2(a)-1\notag\\
=&\nu _2(n(n+2p))-1+\nu _2(a)-0+\nu _2(a)-1\notag\\
=&\nu_2(U_nU_{n+p}U_{n+2p}).
\end{align}
For $q\geq3$, we have $\tau(q)=2$.
Since $q\mid a$ and $\gcd(a,\, b)=1$, we get
\[
\gcd(q,\, \Delta)=\gcd(q,\, a^2+4b)=1.
\]
Now let ${q_1}$ be an odd prime  number dividing $a$ such that
$\nu_{q_1}(a)\geq\nu_{q_1}(n+p)$.
Then $\tau(q_1)=2$ and $q_1\nmid \Delta$.
Note that $U_{\tau(q_1)}=U_2=a$.
By Lemma \ref{lem1}, we get
\begin{align}\label{eq8.14.3}
 &\nu_{q_1}\left(U_{\frac{n(n+p)(n+2p)}{2}\cdot\frac{a}{\gcd(a,\,n+p)}}\right)\notag\\
=&\nu _{q_1}(n(n+p)(n+2p))-\nu _{q_1}(2)+\nu _{q_1}(a)-\nu _{q_1}(\gcd(a,\,n+p))+\nu _{q_1}(U_{\tau(q_1)})\notag\\
=&\nu _{q_1}(n(n+p)(n+2p))-0+\nu _{q_1}(a)-\min\{\nu _{q_1}(a),\, \nu _{q_1}(n+p)\}+\nu _{q_1}(a)\notag\\
=&\nu _{q_1}(n(n+p)(n+2p))+\nu _{q_1}(a)- \nu _{q_1}(n+p)+\nu _{q_1}(a)\notag\\
=&\nu _{q_1}(n(n+2p))+2\nu _{q_1}(U_{\tau(q_1)})\notag\\
=&\nu_{q_1}(U_nU_{n+p}U_{n+2p}).
\end{align}
Let ${q_2}$ be an odd prime  number dividing $a$ such that
$\nu_{q_2}(a)<\nu_{q_2}(n+p)$.
Then $\tau(q_2)=2$ and $q_2\nmid \Delta$.
By Lemma \ref{lem1}, we get
\begin{align}\label{eq8.15.10}
 &\nu_{q_2}\left(U_{\frac{n(n+p)(n+2p)}{2}\cdot\frac{a}{\gcd(a,\,n+p)}}\right)\notag\\
=&\nu _{q_2}(n(n+p)(n+2p))-\nu _{q_2}(2)+\nu _{q_2}(a)-\nu _{q_2}(\gcd(a,\,n+p))+\nu _{q_2}(U_{\tau(q_2)})\notag\\
=&\nu _{q_2}(n(n+p)(n+2p))-0+\nu _{q_2}(a)-\min\{\nu _{q_2}(a),\, \nu _{q_2}(n+p)\}+\nu _{q_2}(a)\notag\\
=&\nu _{q_2}(n(n+2p))+\nu _{q_2}(n+p)+\nu _{q_2}(a)- \nu _{q_2}(a)+\nu _{q_2}(a)\notag\\
>&\nu _{q_2}(n(n+2p))+2\nu _{q_2}(U_{\tau(q_2)})\notag\\
=&\nu_{q_2}(U_nU_{n+p}U_{n+2p}).
\end{align}
From  \eqref{eq5.13.1}-\eqref{eq8.15.10}, for any prime factor $q$ of $a$, we have
\[
\nu_{q}(U_nU_{n+p}U_{n+2p})<\nu_{q}\left(U_{\frac{n(n+p)(n+2p)}{2}\cdot\frac{a}{\gcd(a,\,n+p)}}\right).
\]
Since $\gcd(U_nU_{n+p},U_{n+2p}) =a$, it follows that
\[
\left. U_nU_{n+p}U_{n+2p}\ \middle| \ U_{\frac{n(n+p)(n+2p)}{2}\cdot\frac{a}{\gcd(a,\,n+p)}}\right..
\]
Therefore,
\[
\left. \tau(U_nU_{n+p}U_{n+2p})\middle| \ \frac{n(n+p)(n+2p)}{2}\cdot\frac{a}{\gcd(a,\,n+p)}\right..
\]
We prove by contradiction. Suppose
\[
\tau(U_nU_{n+p}U_{n+2p})\neq\frac{n(n+p)(n+2p)}{2}\cdot\frac{a}{\gcd(a,\,n+p)},
\]
then by \eqref{eq7.23.3}, there exists an integer $t$ such that
\[
\tau(U_nU_{n+p}U_{n+2p})=\frac{n(n+p)(n+2p)}{2}\cdot t,
\]
where
$ t\mid \frac{a}{\gcd(a,\,n+p)}$
and  $t\neq\frac{a}{\gcd(a,\,n+p)}$.
Therefore, there  exists a prime factor $q'$ of $a$ such that
$\nu_{q'}\left(\frac{a}{\gcd(a,\,n+p)}\right)>\nu_{q'}(t)$.
Since
\begin{align}
=&\nu_{q_2}(a)-\nu _{q_2}(\gcd(a,\,n+p))\notag\\
=&\nu_{q_2}(a)-\min\{\nu _{q_2}(a),\, \nu _{q_2}(n+p)\}\notag\\
=&\nu_{q_2}(a)-\nu_{q_2}(a)=0,\notag
\end{align}
it follows that $q'=2$ or $q'=q_1$.
By \eqref{eq8.15.11} and \eqref{eq8.14.3}, we have
\begin{align}
&\nu_2\left(U_{\frac{n(n+p)(n+2p)}{2}\cdot t}\right)\notag\\
=&\nu _2\left(\frac{n(n+p)(n+2p)}{2}\cdot t\right)+\nu _2(a)-1\notag\\
<&\nu _2\left(\frac{n(n+p)(n+2p)}{2}\cdot \frac{a}{\gcd(a,n+p)}\right)+\nu _2(a)-1\notag\\
=&\nu_2(U_nU_{n+p}U_{n+2p}),\notag
\notag\end{align}
and
\begin{align}
&\nu_{q_1}\left(U_{\frac{n(n+p)(n+2p)}{2}\cdot t}\right)\notag\\
=&\nu _{q_1}\left(\frac{n(n+p)(n+2p)}{2}\cdot t\right)+\nu _{q_1}(U_{\tau(q_1)})\notag\\
<&\nu _{q_1}\left(\frac{n(n+p)(n+2p)}{2}\cdot \frac{a}{\gcd(a,n+p)}\right)+\nu _{q_1}(U_{\tau(q_1)})\notag\\
=&\nu_{q_1}(U_nU_{n+p}U_{n+2p}).\notag
\notag\end{align}
This contradicts
$\tau(U_nU_{n+p}U_{n+2p})=\frac{n(n+p)(n+2p)}{2}\cdot t$.
Therefore, we conclude that
\[
\tau(U_nU_{n+p}U_{n+2p})=\frac{n(n+p)(n+2p)}{2}\cdot\frac{a}{\gcd(a,\,n+p)}.
\]

\textbf {Case 3.} $p\mid n$, $p\geq3$ and $2\nmid n$.
For $s\in \{0,\,p,\,2p\}$, it follows from  Lemma \ref{lem5} (iii) that
\[
(n+s)\mid \tau(U_nU_{n+p}U_{n+2p}).
\]
It is obvious that
\[
\gcd(n,n+p)=p,\  \gcd(n,n+2p)=p,\  \gcd(n+p,n+2p)=p.
\]
So there exists an integer $k_1$ such that $\tau(U_nU_{n+p}U_{n+2p})=nk_1$.
This gives $\left.(n+p) \ \middle| \ nk_1 \right.$, yielding $\frac{n+p}{p}\mid k_1$.
Then there exists an integer $k_2$ such that $k_1=\frac{n+p}{p}k_2$,
meaning that $\tau(U_nU_{n+p}U_{n+2p})=\frac{n+p}{p}nk_2$.
Note that  $\left.(n+2p) \ \middle| \  \frac{n+p}{p}nk_2\right.$.
We have
$\left. \frac{n+2p}{p} \ \middle| \ \frac{n+p}{p}\frac{n}{p}k_2\right.$.
It follows that
$\left.\frac{n+2p}{p} \ \middle| \  k_2\right.$.
Therefore, we get
\begin{equation}
\left. \frac{n(n+p)(n+2p)}{p^2} \ \middle| \ \tau(U_nU_{n+p}U_{n+2p}).\right.
\notag \end{equation}
 For $s\in \{0,\,p,\,2p\}$, it is obvious that
 \[
 \left. U_{n+s}\ \middle| \ U_{\frac{n(n+p)(n+2p)}{p^2}U_p^2}\right..
 \]
By Lemma \ref{lem 3} (i), we get that
\[
 \gcd(U_n,\, U_{n+p})=U_p,\  \gcd(U_n,\, U_{n+2p})=U_p,\
\gcd(U_{n+p},\, U_{n+2p})=U_p.
 \]
Let us consider any prime factor $q$ of $U_p$.
Then $\tau(q)\mid p$ and $q\nmid b$.
Since $\tau(q)\geq2$, we have $\tau(q)= p$.
By Lemma \ref{lem2}, we get that $U_p$ is even if and only if
$2\nmid b$, $2\nmid a$, $p=3$.
For $q=2$, by Lemma \ref{lem2}, we obtain
\begin{align}\label{eq7.26.11}
\nu_2\left(U_{\frac{n(n+p)(n+2p)}{p^2}U_p^2}\right)
&=\nu _2(n(n+p)(n+2p))+\nu _2(U_p^2)-\nu_2(p^2)+\nu_2(U_6)-1\notag\\
&=\nu _2(n+p)+\nu _2(U_p^2)-0+\nu_2(U_6)-1\notag\\
&=\nu _2(n+p)+\nu _2(U_3^2)+\nu_2(U_6)-1\notag\\
&=\nu_2(U_nU_{n+p}U_{n+2p}).
 \end{align}
For $q\geq3$, let us consider the following two cases.
If $q\mid \Delta$, then it follows from Lemma \ref{lem 7.24.1} (iii) that  $p=q$.
By Lemma \ref{lem1}, we get
\begin{align}\label{eq7.24.6}
\nu_q\left(U_{\frac{n(n+p)(n+2p)}{p^2}U_p^2}\right)
&=\nu _q(n(n+p)(n+2p))+\nu _q(U_p^2)-\nu_q(p^2)+\nu_q(U_q)-1\notag\\
&=\nu _q(n(n+p)(n+2p))+3\nu _q(U_{q})-3\notag\\
&=\nu_q(U_nU_{n+p}U_{n+2p}).
 \end{align}
Next we can assume that
$q\nmid \Delta$. By Lemma \ref{lem 7.24.1} (iii), we get $p\neq q$.
By Lemma \ref{lem1}, we have
\begin{align}\label{eq7.24.7}
\nu_q\left(U_{\frac{n(n+p)(n+2p)}{p^2}U_p^2}\right)
&=\nu _q(n(n+p)(n+2p))+\nu _q(U_p^2)-\nu_q(p^2)+\nu_q(U_{\tau(q)})\notag\\
&=\nu _q(n(n+p)(n+2p))+3\nu _q\left(U_{\tau(q)}\right)\notag\\
&=\nu_q(U_nU_{n+p}U_{n+2p}) .
\end{align}
Together with \eqref{eq7.26.11}--\eqref{eq7.24.7},
we have
\[
\left. U_nU_{n+p}U_{n+2p} \ \middle| \  U_{\frac{n(n+p)(n+2p)}{p^2}U_p^2}.\right.
\]
This gives
\[
\left. \tau(U_nU_{n+p}U_{n+2p}) \ \middle| \ \frac{n(n+p)(n+2p)}{p^2}U_p^2.\right.
\]
Now, we claim that
\[
\tau(U_nU_{n+p}U_{n+2p})=\frac{n(n+p)(n+2p)}{p^2}U_p^2.
\]
Let us assume the opposite, namely that there exists an integer $t$ such that $\tau(U_nU_{n+p}U_{n+2p})=\frac{n(n+p)(n+2p)}{p^2}t $, where $t\mid U_p^2$ but $t\neq U_p^2$.
It follows that there exists a prime number $q_1$ such that $\nu_{q_1}(U_p^2)> \nu_{q_1}(t)$.
By \eqref{eq7.26.11}--\eqref{eq7.24.7}, we know that
\[
\nu_{q_1}\left(U_{\frac{n(n+p)(n+2p)}{p^2}t}\right)< \nu_{q_1}\left(U_{\frac{n(n+p)(n+2p)}{p^2}U_p^2}\right)=\nu_{q_1}\left(U_nU_{n+p}U_{n+2p}\right),
\]
which is a contradiction. The claim follows.

\textbf {Case 4.} $p\mid n$, $p\geq3$ and $2\mid n$.
For $s\in \{0,\,p,\,2p\}$, it follows from  Lemma \ref{lem5} (iii) that
\[
(n+s)\mid \tau(U_nU_{n+p}U_{n+2p}).
\]
It is obvious that
\[
\gcd(n,n+p)=p,\  \gcd(n,n+2p)=2p,\  \gcd(n+p,n+2p)=p.
\]
So there exists an integer $k_3$ such that $\tau(U_nU_{n+p}U_{n+2p})=nk_3$.
This gives $\left.(n+p) \ \middle| \ nk_3\right.$, yielding $\frac{n+p}{p}\mid k_3$.
Then there exists an integer $k_4$ such that $k_3=\frac{n+p}{p}k_4$,
meaning that $\tau(U_nU_{n+p}U_{n+2p})=\frac{n+p}{p}nk_4$.
Since $\left.\frac{n+2p}{2p}\ \middle| \ \frac{n+p}{p}\frac{n}{2p}k_4\right.$, we get
$\left.\frac{n+2p}{2p} \ \middle| \  k_4\right.$.
Therefore, we get
\begin{equation}\label{eq7.24.8}
\left. \frac{n(n+p)(n+2p)}{2p^2} \ \middle| \ \tau(U_nU_{n+p}U_{n+2p}).\right.
\end{equation}
For $s\in \{0,\,p,\,2p\}$, it is obvious that
 \[
 \left. U_{n+s}\ \middle| \ U_{\frac{n(n+p)(n+2p)}{2p^2}\cdot\frac{U_p^2V_p}{\gcd\left(V_p,\, (n+p)/p\right)}}.\right.
 \]
By Lemma \ref{lem 3} (i), we get
\[
 \gcd(U_n,\, U_{n+p})=U_p,\, \gcd(U_n,\, U_{n+2p})=U_{2p}=U_pV_p,\, \gcd(U_{n+p},\, U_{n+2p})=U_p.
 \]
 By Lemma \ref{lem 3} (iii), we know that $\gcd(U_p,\, V_p)=1$ or $2$.
 If $2\nmid a$, $2\nmid b$, $p=3$ or $2\mid a$, then at least one of the numbers $U_p$ and $V_p$ is even.
Otherwise,  we find that $U_p$ and $V_p$ both are odd.
For $2\mid a$, it follows from  Lemmas \ref{lem2} and  \ref{lemma:3.311} that
\begin{align}\label{eq7.25.3}
&\nu_{2}\left(U_{\frac{n(n+p)(n+2p)}{2p^2}\cdot\frac{U_p^2V_p}{\gcd(V_p,\, (n+p)/p)}}\right)\notag\\
=&\nu _2(n(n+p)(n+2p))-\nu_2(2p^2)+\nu_2(U_p^2V_p)-\nu_2(\gcd(V_p,\, (n+p)/p))+\nu_2(a)-1\notag\\
=&\nu _2(n(n+2p))-1+\nu_2(V_p)-0+\nu_2(a)-1\notag\\
=&\nu _2(n(n+2p))-1+\nu_2(a)+\nu_2(a)-1\notag\\
=&\nu_2(U_nU_{n+p}U_{n+2p}).
 \end{align}
For $2\nmid a$, $2\nmid b$ and $p=3$, we obtain
\begin{align}\label{eq7.26.12}
&\nu_{2}\left(U_{\frac{n(n+p)(n+2p)}{2p^2}\cdot\frac{U_p^2V_p}{\gcd(V_p,\, (n+p)/p)}}\right)\notag\\
=&\nu _2(n(n+p)(n+2p))+\nu _2(U_p^2V_p)-\nu_2(2p^2)-\nu_2(\gcd(V_p,\, (n+p)/p))+\nu_2(U_6)-1\notag\\
=&\nu _2(n(n+2p))+\nu _2(U_3U_{6})+\nu_2(U_6)-2\notag\\
=&\nu _2(n(n+2p))+2\nu_2(U_6)-2+\nu _2(U_3)\notag\\
=&\nu_2(U_nU_{n+p}U_{n+2p}).
 \end{align}
If $U_p$ and $V_p$ both are odd, then $\gcd(U_p,\, V_p)=1$.
Let $q$ be an odd prime number  dividing $U_p$. Then $q\nmid V_p$.
If $q\mid \Delta$, then $\tau(q)=p$ and $p=q$.
By Lemma \ref{lem1}, we get
\begin{align}\label{eq7.25.1}
&\nu_q\left(U_{\frac{n(n+p)(n+2p)}{2p^2}\cdot\frac{U_p^2V_p}{\gcd(V_p,\, (n+p)/p)}}\right)\notag\\
=&\nu _q(n(n+p)(n+2p))-\nu_q(2p^2)+\nu_q(U_p^2V_p)-\nu_q(\gcd(V_p,\, (n+p)/p))+\nu _q(U_{q})-1\notag\\
=&\nu _q(n(n+p)(n+2p))+3\nu _q(U_{q})-3\notag\\
=&\nu_q(U_nU_{n+p}U_{n+2p}).
 \end{align}
If $q\nmid \Delta$, then $\tau(q)=p$ and $p\neq q$.
By Lemma \ref{lem1}, we get
\begin{align}\label{eq7.25.2}
&\nu_q\left(U_{\frac{n(n+p)(n+2p)}{2p^2}\cdot\frac{U_p^2V_p}{\gcd(V_p,\, (n+p)/p)}}\right)\notag\\
=&\nu _q(n(n+p)(n+2p))-\nu_q(2p^2)+\nu_q(U_p^2V_p)-\nu_q(\gcd(V_p,\, (n+p)/p))+\nu _q(U_{\tau(q)})\notag\\
=&\nu _q(n(n+p)(n+2p))-0+\nu_q(U_p^2)+\nu_q(U_p)\notag\\
=&\nu _q(n(n+p)(n+2p))+3\nu _q(U_{\tau(q)})\notag\\
=&\nu_q(U_nU_{n+p}U_{n+2p}).
 \end{align}
Now let ${q_1}$ be an odd prime  number dividing $V_p$ such that
$\nu_{q_1}\left((n+p)/p\right)\leq \nu_{q_1}(V_p)$.
Then $\nu_{q_1}(V_p)>0$ and ${q_1}\nmid U_p$.
By Lemma \ref{lemma:3.312}, we get ${q_1}\nmid \Delta$, $\tau(q_1)\nmid p$ and $\tau(q_1)\mid 2p$.
Since $2p\mid n$ and $2p\mid (n+2p)$, it follows that $\tau(q_1)\mid n$, $\tau(q_1)\mid (n+2p)$ and $\tau(q_1)\nmid (n+p)$.
By Lemmas \ref{lem1} and  \ref{lemma:3.312}, we obtain
\begin{align}\label{eq7.25.4}
&\nu_{q_1}\left(U_{\frac{n(n+p)(n+2p)}{2p^2}\cdot\frac{U_p^2V_p}{\gcd(V_p,\, (n+p)/p)}}\right)\notag\\
=&\nu _{q_1}(n(n+p)(n+2p))-\nu_{q_1}(2p^2)+\nu_{q_1}(U_p^2V_p)-\nu_{q_1}(\gcd(V_p,\, (n+p)/p))+\nu _{q_1}(U_{\tau({q_1})})\notag\\
=&\nu _{q_1}(n(n+p)(n+2p))-2\nu_{q_1}(p)+\nu_{q_1}(V_p)-\min\{\nu_{q_1}(V_p),\nu_{q_1}\left((n+p)/p\right)\}+\nu _{q_1}(U_{\tau({q_1})})\notag\\
=&\nu _{q_1}(n(n+2p))+\nu_{q_1}(n+p)-2\nu_{q_1}(p)+(\nu_{q_1}(p)+\nu _{q_1}(U_{\tau({q_1})}))-\nu_{q_1}((n+p)/p)+\nu _{q_1}(U_{\tau({q_1})})\notag\\
=&\nu _{q_1}(n(n+2p))+2\nu _{q_1}(U_{\tau({q_1})})\notag\\
=&\nu_{q_1}(U_nU_{n+p}U_{n+2p}).
 \end{align}
Let ${q_2}$ be an odd prime  number dividing $V_p$ such that
$\nu_{q_2}\left((n+p)/p\right)> \nu_{q_2}(V_p)$.
Then $\nu_{q_2}(V_p)>0$ and ${q_2}\nmid U_p$.
Similarly, we have
\[
{q_2}\nmid \Delta,\ \tau(q_2)\mid n,\  \tau(q_2)\mid (n+2p),\  \tau(q_2)\nmid (n+p).
\]
So, we obtain
\begin{align}\label{eq8.15.2}
&\nu_{q_2}\left(U_{\frac{n(n+p)(n+2p)}{2p^2}\cdot\frac{U_p^2V_p}{\gcd\left(V_p,\, (n+p)/p\right)}}\right)\notag\\
=&\nu _{q_2}(n(n+p)(n+2p))-\nu_{q_2}(2p^2)+\nu_{q_2}(U_p^2V_p)-\nu_{q_2}\left(\gcd\left(V_p,\, (n+p)/p\right)\right)+\nu _{q_2}(U_{\tau({q_2})})\notag\\
=&\nu _{q_2}(n(n+p)(n+2p))-2\nu_{q_2}(p)+\nu_{q_2}(V_p)-\min\{\nu_{q_2}(V_p),\nu_{q_2}\left((n+p)/p\right)\}+\nu _{q_2}(U_{\tau({q_2})})\notag\\
=&\nu _{q_2}(n(n+2p))+\nu_{q_2}(n+p)-2\nu_{q_2}(p)+\nu_{q_2}(V_p)-\nu_{q_2}(V_p)+\nu _{q_2}(U_{\tau({q_2})})\notag\\
> &\nu _{q_2}(n(n+2p))+\nu_{q_2}(V_p)-\nu_{q_2}(p)+\nu _{q_2}(U_{\tau({q_2})})\notag\\
=&\nu _{q_2}(n(n+2p))+2\nu _{q_2}(U_{\tau({q_2})})\notag\\
=&\nu_{q_2}(U_nU_{n+p}U_{n+2p}).
 \end{align}

By \eqref{eq7.25.3}-\eqref{eq8.15.2}, we have
\[
\left.U_nU_{n+p}U_{n+2p}\middle|\ U_{\frac{n(n+p)(n+2p)}{2p^2}\cdot\frac{U_p^2V_p}{\gcd(V_p,(n+p)/p)}}\right..
\]
So,
\[
\left. \tau(U_nU_{n+p}U_{n+2p})\middle|\  \frac{n(n+p)(n+2p)}{2p^2}\cdot\frac{U_p^2V_p}{\gcd(V_p,\, (n+p)/p)}\right..
\]
By \eqref{eq7.24.8},  there exists an integer $t$ such that
$ \tau(U_nU_{n+p}U_{n+2p})=\frac{n(n+p)(n+2p)}{2p^2}t$,
where
$t\mid \frac{U_p^2V_p}{\gcd(V_p,\, (n+p)/p)}$.
We prove by contradiction.
Suppose $t\neq \frac{U_p^2V_p}{\gcd(V_p,\, (n+p)/p)}$, there exists a prime factor $q'$ such that
\[
\nu_{q'}\left(\frac{U_p^2V_p}{\gcd(V_p,\, (n+p)/p)}\right)>\nu_{q'}(t),
\]
where $q'$ is a prime factor of $U_p$ or $V_p$.
By \eqref{eq7.25.3}-\eqref{eq7.25.4}, we have
$q'\neq 2$, $q'\neq q$ and  $q'\neq q_1$.
So, $q'= q_2$. Since
 \begin{align}
&\nu_{q_2}\left(\frac{U_p^2V_p}{\gcd\left(V_p,\, (n+p)/p\right)}\right)\notag\\
=&\nu _{q_2}(U_p^2V_p)-\nu_{q_2}\left(\gcd\left(V_p,\, (n+p)/p\right)\right)\notag\\
=&\nu _{q_2}(V_p)-\min\{\nu_{q_2}(V_p),\nu_{q_2}\left((n+p)/p\right)\}\notag\\
=&\nu _{q_2}(V_p)-\nu _{q_2}(V_p)=0,\notag
\notag \end{align}
this contradicts  $\nu_{q'}\left(\frac{U_p^2V_p}{\gcd(V_p,\, (n+p)/p)}\right)>\nu_{q'}(t)$.
Therefore, there does not exist a prime number $q'$ satisfying the above conditions.
So, $t= \frac{U_p^2V_p}{\gcd(V_p,\, (n+p)/p)}$.
Therefore, we conclude that
\[
\tau(U_nU_{n+p}U_{n+2p})=\frac{n(n+p)(n+2p)}{2p^2}\cdot\frac{U_p^2V_p}{\gcd(V_p,\, (n+p)/p)}.
\]

\end{prof}


\begin{thebibliography}{}
\bibitem{Irmak} N. Irmak, P. K. Ray,  The order of appearance of the product of two Fibonacci and Lucas numbers, \emph{Acta Math. Hungar.}, \textbf{162}  (2020), 527--538.
\bibitem{Khaochim1} N. Khaochim, P. Pongsriiam, On the order of appearance of products of Fibonacci numbers, \emph{Contrib. Discrete Math.}, \textbf{13}  (2018), 45--62.
\bibitem{Khaochim2} N. Khaochim, P. Pongsriiam, The general case on the order of appearance of product of consecutive Lucas numbers, \emph{Acta Math. Univ. Comenian. (N.S.)}, \textbf{87}  (2018), 277--289.
\bibitem{5Marques}D. Marques, The order of appearance of integers at most one away from Fibonacci numbers, \emph{Fibonacci Quart.}, \textbf{50} (2012), 36--43.
\bibitem{6Marques} D. Marques, The order of appearance of product of consecutive Fibonacci numbers, \emph{Fibonacci Quart.}, \textbf{50}  (2012), 132--139.
\bibitem{7Marques}D. Marques, The order of appearance of powers of Fibonacci and Lucas numbers, \emph{Fibonacci Quart.}, \textbf{50}  (2012), 239--245.
\bibitem{8Marques} D. Marques, The order of appearance of the product of consecutive Lucas numbers, \emph{Fibonacci Quart.}, \textbf{51}  (2013), 38--43.
\bibitem{McDaniel} W. L. McDaniel, The g.c.d in Lucas sequences and Lehmer number sequences, \emph{Fibonacci Quart.}, \textbf{29}  (1991), 24--29.
\bibitem{Onphaeng} K. Onphaeng, P. Pongsriiam, Exact divisibility by powers of the integers in the Lucas sequences of the first and second kinds, \emph{AIMS Math.}, \textbf{6}  (2021), 11733--11748.
\bibitem{Pongsriiam} P. Pongsriiam, A complete formula for the order of appearance of the powers of Lucas numbers, \emph{Commun. Korean Math. Soc. }, \textbf{31} (2016), 447--450.
\bibitem{Ray} P. K. Ray,  N. Irmak, B. K. Patel,  The rank of apparition of powers of Lucas sequence, \emph{Turkish J. Math.}, \textbf{42} (2018), 1566--1570.
\bibitem{Renault} M. Renault, The period, rank, and order of the $(a,b)$-Fibonacci sequence mod $m$, \emph{Math. Mag.}, \textbf{86} (2013), 372--380.
\bibitem{Ribenboim} P. Ribenboim,  My numbers, my friends, Springer-Verlag, New York, 2000.
\bibitem{Sanna}  C. Sanna,   The $p$-Adic valuation of Lucas sequences, \emph{Fibonacci Quart.}, \textbf{54}  (2016), 118--124.
\bibitem{2Trojovsky} P. Trojovsk\'{y}, On the order of appearance of the difference of two Lucas numbers, \emph{Miskolc Math. Notes}, \textbf{19}  (2018), 641--648.
\bibitem{1Trojovsky} P. Trojovsk\'{y},  The order of appearance of the sum and difference between two Fibonacci numbers, \emph{Asian-Eur. J. Math.}, \textbf{12}  (2019), no. 3, 1950046, 10 pp.

\end{thebibliography}
\end{document}